\documentclass[10pt]{article}

\usepackage{amsmath}
\usepackage{amssymb}
\usepackage{amsthm}
\usepackage{bbm}
\usepackage{graphicx}
\usepackage{color}
\usepackage{epsfig}

\newtheorem{theorem}{Theorem}
\newtheorem{prop}[theorem]{Proposition}
\newtheorem{lemma}[theorem]{Lemma}
\newtheorem{cor}[theorem]{Corollary}

\newtheorem{question}{Question}

\newtheorem{conj}[question]{Conjecture}

\newcommand{\1}{{\mathbbmss 1}}

\newcommand{\cI}{{\mathcal I}}
\newcommand{\cL}{{\mathcal L}}

\DeclareMathOperator{\disc}{disc}

\DeclareMathOperator{\run}{run}
\DeclareMathOperator{\rp}{rp}
\DeclareMathOperator{\sk}{sk}

\allowdisplaybreaks

\title{The Discrepancy of the Lex-Least De Bruijn Sequence}
\author{Joshua Cooper and Christine Heitsch}

\begin{document}

\maketitle

\begin{abstract} We answer the following question of R.~L.~Graham: What is the discrepancy of the lexicographically-least binary de Bruijn sequence?  Here, ``discrepancy'' refers to the maximum (absolute) difference between the number of ones and the number of zeros in any initial segment of the sequence.  We show that the answer is $\Theta(2^n \log n/n)$.
\end{abstract}

\section{Introduction}

A {\it binary de Bruijn sequence of order $k$} is a word $a_1 \cdots a_{2^k}$ over the alphabet $\{0,1\}$ that contains every $k$-word exactly once as a subword when
the indices are interpreted cyclically.  It is well known (see, e.g., \cite{L01}) that the number of de Bruijn cycles of order $k$ is given by
$$
2^{2^{k-1}-k}.
$$
Among these is the ``Ford sequence''\footnote{See the excellent survey \cite{F82} for a history of this and related sequences.  The eponym, due to Fredricksen, refers to a 1957 unpublished manuscript of Ford (\cite{F57}).  However, subsequent research has revealed earlier references.  In \cite{F82}, the author proposes that a 1934 paper of Martin (\cite{M34}) is the earliest appearance.}, the remarkable cyclic binary word which is
\begin{enumerate}
\item the lexicographic least de Bruijn sequence,
\item the result of applying the least-first greedy algorithm to constructing a de Bruijn sequence (starting with $1^k$),
\item the result of concatenating all ``Lyndon'' words (lexicographically minimal representatives of free conjugacy classes) of each length dividing $k$ in lexicographic order, and
\item the de Bruijn sequence generated by a shift register whose truth table has minimum weight.
\end{enumerate}
Since the greedy algorithm uses 0's before 1's whenever possible, it is natural to suspect that this special sequence has an excess of $0$'s early on, i.e, the difference between the number of $0$'s and $1$'s in initial segments is large.  Indeed, Huang comments in \cite{H90} that
\begin{quote}
The ``prefer one'' algorithm proposed by Fredricksen joins the pure cycles of [a] circulating register (CR) in order according to the weights of the $n$-tuples... so some part of the sequence may contain many heavily weighted $n$-tuples and it leads to a bad local 0-1 balance.
\end{quote}
R.~L.~Graham therefore asks for the maximum ``discrepancy.''  In the present note, we show that it has order $2^n \log n/n$.\\

Define the equivalence relation $\sim$ (``conjugacy'') on binary words by setting $xy \sim yx$ for any $x,y \in \{0,1\}^\ast$.  For a word $w \in \{0,1\}^\ast$, define $w^\circ$ to be the lexicographic least element of the $\sim$-equivalence class $[\![w]\!]$ of $w$.  If $w$ is aperiodic (i.e., if $w = xy$ with $x,y \neq \epsilon$, then $w \neq yx$), then $w^\circ$ is called a ``Lyndon word.''  Then the lexicographically least binary order-$n$ de Bruijn sequence $\cL_n$ consists of the concatenation of all Lyndon words of length dividing $n$, {\it in lexicographic order}.

For a word $w \in \{0,1\}^\ast$, write $w_k$ for its $k^\textrm{th}$ symbol from left to right, starting with zero.  Then we define the {\it discrepancy} of $w$ to be
$$
\disc(w) = \max_M \left |\sum_{k=0}^M (-1)^{w_k} \right |.
$$

\begin{theorem} $\disc(\cL_n) = \Theta(2^n \log n / n)$. \label{thm:main}
\end{theorem}

We conjecture a slightly stronger statement:

\begin{conj} There is some $C$ so that $\lim_{n\rightarrow \infty} \frac{n \disc(\cL_n)}{2^n \log n} = C$.
\end{conj}

Our argument will estimate the discrepancy of $\cL_n$ by considering
substrings consisting of Lyndon words $w^{\circ}$ grouped by the
length $k$ of their $0^{k}1$ prefix.
For $0 < k < n$, let $S_k$ be the set of binary words of length $n$
containing the subword $0^k$ but not the subword $0^{k+1}$.
Then the elements of $S_k$ are precisely those $w$ so that $w^\circ$
begins with $0^k$.
Define $S_k^\circ = \{w^\circ : w \in S_k\}$, and let $\ell_k$ be the
concatenation of the elements of $S_k^\circ$ in lexicographic order.
Since the elements of $S_k^\circ$ precede those of $S_{k-1}^\circ$
in the lexicographic order, this means that
$$
\cL_n = 0 \cdot \left ( \prod_{k=1}^{n-1} \ell_{n-k} \right ) \cdot 1,
$$
as long as $n$ is prime.

For a binary string $w$ of length $n$, we define the {\it skew} of $w$
to be
\[ \sk(w) = \sum_{i = 0}^{n-1} (-1)^{w_{i}} \]
so that
$$
\disc(\cL_n) = \max_{1 \leq t \leq n-2} \left ( 1 + \sum_{k=1}^t \sk(\ell_{n-k}) + \disc(\ell_{n-t-1}) \right )
$$
when $n$ is prime.  This will allow us to bound the discrepancy of $\cL_n$.

\section{Preliminaries}

Define $\alpha_k(n)$ to be the number of elements of $\{0,1\}^n$ containing no subword $0^k$, and let $\beta_k(n)$ be defined by
$$
\beta_k(n) = \sum_{\substack{w \in \{0,1\}^n \\ 0^k \not \in w}} \sk(w).
$$
For the remainder of this section, we fix a $k \geq 2$.

\begin{lemma} \label{lem:recurrence} The sequences $a_n = \alpha_k(n)$ and $b_n = \beta_k(n)$ satisfy:
\begin{enumerate}
\item $a_n = \sum_{j=1}^k a_{n-j}$ for $n \geq k$, and
\item $b_n = \sum_{j=1}^k [(j-2) a_{n-j} + b_{n-j}]$ for $n \geq k$.
\end{enumerate}
Furthermore, $a_j = 2^j$ for $0 \leq j < k$ and $b_j = 0$ for $0 \leq j < k$.
\end{lemma}
\begin{proof} Both recurrences follow from the following consideration: any string of length at least $k$ not containing a subword $0^k$ has a left-most $1$.  Therefore, we may partition the $0^k$-free sequences into those which begin with a string of the form $0^j 1$ for $0 \leq j < k$.  The ``base case'' formulas trivially follow from the fact that {\it every} string of length less than $k$ is $0^k$-free.
\end{proof}

\begin{lemma} \label{lem:recurrencebound} For $n-1 \geq k \geq 3$,
$$
a_{n-1} = k + \sum_{j=3}^k (j-2) a_{n-j} + (k-1) \sum_{j=0}^{n-k-1} a_j.
$$
\end{lemma}
\begin{proof} We proceed by induction.  First, we verify that $a_{k} = k + \sum_{j=3}^{k} (j-2) a_{k+1-j} + (k-1)a_0$.  Note that, by the ``base case'' part of Lemma \ref{lem:recurrence}, $a_j = 2^j$ in the relevant range, except that $a_k = 2^k-1$.  Therefore,
\begin{align*}
k + \sum_{j=3}^k (j-2) a_{k+1-j} + (k-1) a_0 &= k + \sum_{j=3}^k (j-2) 2^{k+1-j} + k-1\\
&= \sum_{j=1}^{k-2} j 2^{k-j-1}  + 2k-1\\
&= 2^{k-2} \sum_{j=1}^{k-2} j 2^{-(j-1)} + 2k-1  \\
&= 2^{k-2} (4 - k 2^{-k+3} )  + 2k-1\\
&= 2^{k} - 2k  + 2k-1 \\
& = 2^{k} - 1 = a_{k}.
\end{align*}

Now, suppose the statement holds for $n$.  Applying the first recurrence in Lemma \ref{lem:recurrence},
\begin{align*}
a_n &= \sum_{j=1}^k a_{n-j} \\
&= a_{n-1} + \sum_{j=2}^k a_{n-j} \\
&= k + \sum_{j=3}^k (j-2) a_{n-j} + (k-1) \sum_{j=0}^{n-k-1} a_j + \sum_{j=2}^k a_{n-j} \\
&= k + \sum_{j=2}^k (j-1) a_{n-j} + (k-1) \sum_{j=0}^{n-k-1} a_j \\
&= k + \sum_{j=3}^{k} (j-2) a_{n+1-j} + (k-1) a_{n-k} + (k-1) \sum_{j=0}^{n-k-1} a_j \\
&= k + \sum_{j=3}^{k} (j-2) a_{n+1-j} + (k-1) \sum_{j=0}^{n-k} a_j.
\end{align*}
\end{proof}

\begin{cor} \label{cor:justnegative} $b_n < 0$ for all $n-1 \geq k \geq 3$.
\end{cor}
\begin{proof} If we combine the recurrence for $b_n$ from Lemma \ref{lem:recurrence} with the above Lemma \ref{lem:recurrencebound},
\begin{align}
b_n &= \sum_{j=1}^k [ (j-2) a_{n-j} + b_{n-j} ] \nonumber \\
&= -a_{n-1} + \sum_{j=3}^k (j-2) a_{n-j} + \sum_{j=1}^k b_{n-j} \nonumber \\
&= -k - (k-1) \sum_{j=0}^{n-k-1} a_j + \sum_{j=1}^k b_{n-j} < 0 \label{eq1},
\end{align}
by induction.
\end{proof}

Let $\rho_k$ be the largest (in absolute value) root of the polynomial $g(z) = z^{k+1} - 2z^k + 1$.  It is proven in \cite{O95} that $\rho_k$ is real, lies between $5/3$ and $2$, and is unique in these respects.  It is also shown in \cite{O95} that $\rho_k \rightarrow 2$ as $k \rightarrow \infty$.  Note that
$$
z^k - \sum_{j=0}^{k-1} z^j = \frac{z^{k+1} - 2z^k + 1}{z-1},
$$
so that $\rho_k$ is a root of the left-hand polynomial $f(z)$ here as well.  Since $f(z)$ is the characteristic polynomial for the recurrence that the $a_n$ satisfy, $\rho_k$ is the growth rate of the $a_n$, i.e., $\lim_{n \rightarrow \infty} \log a_n / n = \rho_k$.

\begin{lemma} \label{lem:lowerbound} For all $n \geq 1$, $a_{n} \geq \rho_k a_{n-1}$.
\end{lemma}
\begin{proof} Since $\rho_k < 2$, and $a_n = 2^n$ for $0 \leq n < k$, the claimed bound holds for $n$ in this range.  Suppose it holds for all $n < N$. Then by Lemma \ref{lem:recurrence},
\begin{align*}
a_n &= \sum_{j=1}^k a_{n-j} \\
&\geq \sum_{j=1}^k \rho_k a_{n-j-1} \\
&= \rho_k a_{n-1}.
\end{align*}
\end{proof}

\begin{lemma} \label{lem:finalbound} For $k \geq 4$ and all $n \geq k$, $b_n \geq -2 k a_n/3$.
\end{lemma}
\begin{proof} By (\ref{eq1}),
\begin{align*}
b_n &=  -k - (k-1) \sum_{j=0}^{n-k-1} a_j + \sum_{j=1}^k b_{n-j}.
\end{align*}
If we suppose that $b_{j} \geq -\gamma k a_j$ for all $j<n$, then
\begin{align*}
b_n &\geq  -k - k \sum_{j=0}^{n-k-1} a_j - \gamma k \sum_{j=1}^k a_{n-j} \\
&=  -k - k \sum_{j=0}^{n-k-1} a_j - \gamma k a_n.
\end{align*}
By iterating Lemma \ref{lem:lowerbound}, we have
\begin{align*}
b_n &\geq -k - k \sum_{j=0}^{n-k-1} \rho_k^{j-n} a_n - \gamma k a_n \\
&\geq - a_n k \left ( \frac{1}{a_n} + \sum_{j=0}^{\infty} \rho_k^{j-n} + \gamma \right ) \\
&= -a_n k \left (\frac{1}{a_n} + \frac{\rho_k^{-n}}{1 - \rho_k^{-1}} + \gamma \right )\\
&\geq -a_n k \left (\frac{1}{a_n} + \frac{5}{2} \rho_k^{-n} + \gamma \right )
\end{align*}
We may begin by taking $\gamma = \frac{2k-1}{k (2^{k+1}-3)} \leq \frac{7}{116}$ by considering $a_{k+1} = 2^{k+1}-3$ and $b_{k+1} = 1-2k$.  Then, $\gamma$ increases by at most
\begin{align*}
\sum_{n = k+1}^\infty \left (\frac{1}{a_n} + \frac{5}{2} \rho_k^{-n} \right ) & \leq \sum_{n = k+1}^\infty \frac{1}{\rho_k^{n-k} a_k} + \frac{5}{2} \sum_{n = k+1}^\infty \rho_k^{-n} \\
& = \frac{\rho_k^k}{2^k-1} \sum_{n = k+1}^\infty \rho_k^{-n} + \frac{5}{2} \sum_{n = k+1}^\infty \rho_k^{-n} \\
& = \left ( \frac{\rho_k^k}{2^k-1} + \frac{5}{2} \right ) \rho_k^{-k-1} \sum_{n = 0}^\infty \rho_k^{-n} \\
& = \left ( \frac{\rho_k^{-1}}{2^k-1} + \frac{5}{2\rho_k^{k+1} } \right ) \cdot \frac{1}{1-\rho_k^{-1}} \\
& \leq \left ( \frac{3}{5 \cdot 15} + \frac{5}{2 (5/3)^{5} } \right ) \cdot \frac{5}{2} = \frac{293}{500}.
\end{align*}
The conclusion follows for all $n \geq k+1$, since $\frac{293}{500} + \frac{7}{116} = \frac{2343}{3625} \leq \frac{2}{3}$.  It is also easy to verify that $b_k \geq - 2 k a_k/3$.
\end{proof}

\section{Main Result}

Here we prove Theorem \ref{thm:main} stated in the introduction.

\begin{prop} \label{ineq}
For $4 \leq k < n$ and $n$ prime,
$$
\frac{k}{3} - 2 \leq \frac{\sk(\ell_k)}{\alpha_{k+1}(n-k-2)} \leq 2k-3.
$$
\end{prop}

\begin{proof}

The set $S_k$ contains each sequence of the form $0^{k}1w$ where $w$ is a $0^k$-free word of length $n-k-1$.  However, the quantity $\sk(S_k)$ is not quite the sum of the skews of all $0^k$-free sequences of length $n-k-1$ prefixed by $0^k1$: it must include all elements of $S_k^\circ$, not just those that have prefix $0^k$ and contain no other runs $0^k$.  For each word $w$ of length $n$ which contains more than one run of the form $0^k$, but no runs of the form $0^{k+1}$, only one of its conjugates (namely, $w^\circ$) appears in $S_k^\circ$.  Define $\run(w)$ to be the maximum $k$ so that $0^k \in w$, and let $\rho_k(w)$ be the number of subwords of the form $0^k$ in $w$, where $\run(w) = k$.  (Set $\rho_k(w) = 0$ otherwise.)  Since we may assume that each $w$ is aperiodic, this means that
\begin{align*}
\sk(\ell_k) &= \sum_{w \in S_k^\circ} \sk(w) \\
 &= \sum_{\substack{w \in \{0,1\}^n \\ \run(w) = k}} \1(w = w^\circ) \sk(w) \\
 &= \sum_{t \geq 0} \sum_{\substack{w \in \{0,1\}^{n-k-2} \\ \rho_k(w) = t}} \frac{\sk(0^k1w1)}{t+1} \\
 &= \sum_{t \geq 0} \frac{1}{t+1} \sum_{\substack{w \in \{0,1\}^{n-k-2} \\ \rho_k(w) = t}} (k - 2 + \sk(w)).
\end{align*}
Define the ``run-print'' $\rp(w)$ of a word $w \in \{0,1\}$ with $\run(w) = k$ to be the set of indices $j \in [n]$ so that $w$ has a run $0^k$ starting at index $j$.  Then we may write
\begin{align*}
\sk(\ell_k) &= \sum_{t \geq 0} \frac{1}{t+1} \sum_{\substack{w \in \{0,1\}^{n-k-2} \\ \rho_k(w) = t}} (k - 2 + \sk(w)) \\
&= \sum_{t \geq 0} \frac{1}{t+1} \sum_{\substack{w \in \{0,1\}^{n-k-2} \\ \rho_k(w) = t}} (k - 2) \\
&\qquad + \sum_{t \geq 0} \frac{1}{t+1} \sum_{S \in \binom{[n-k-2]}{t}} \sum_{\substack{w \in \{0,1\}^{n-k-2} \\ \rp(w) = S}} \sk(w).
\end{align*}
Now, for a given $S$ of cardinality $t$ and $w$ with $\rp(w)=S$, there is a $0^k$ run starting at location $s$ for each $s \in S$.  Each such run is bounded on both sides by a $1$.  In between the runs are intervals, the sum over whose skews is nonpositive, by Corollary \ref{cor:justnegative}.  Therefore,
\begin{align*}
\sum_{\substack{w \in \{0,1\}^{n-k-2} \\ \rp(w) = S}} \sk(w) &\leq \sum_{\substack{w \in \{0,1\}^{n-k-2} \\ \rp(w) = S}} t(k-1) ,
\end{align*}
so we have
\begin{align*}
\sk(\ell_k) &\leq \sum_{t \geq 0} \frac{1}{t+1} \sum_{\substack{w \in \{0,1\}^{n-k-2} \\ \rho_k(w) = t}} (k - 2) \\
&\qquad + \sum_{t \geq 0} \frac{1}{t+1} \sum_{S \in \binom{[n-k-2]}{t}} \sum_{\substack{w \in \{0,1\}^{n-k-2} \\ \rp(w) = S}} t(k-1) \\
&< \sum_{t \geq 0} \sum_{\substack{w \in \{0,1\}^{n-k-2} \\ \rho_k(w) = t}} (k - 2) + (k-1) \sum_{t \geq 0} \sum_{S \in \binom{[n-k-2]}{t}} \sum_{\substack{w \in \{0,1\}^{n-k-2} \\ \rp(w) = S}} 1 \\
&= (k-2) \alpha_{k+1}(n-k-2) + (k-1) \alpha_{k+1}(n-k-2) \\
&= (2k-3) \alpha_{k+1}(n-k-2).
\end{align*}
On the other hand, by Lemma \ref{lem:finalbound},
\begin{align*}
\sk(\ell_k) &= \sum_{t \geq 0} \frac{1}{t+1} \sum_{\substack{w \in \{0,1\}^{n-k-2} \\ \rho_k(w) = t}} (k - 2 + \sk(w)) \\
&\geq \sum_{\substack{w \in \{0,1\}^{n-k-2} \\ \rho_k(w) = 0}} (k - 2) + \sum_{\substack{w \in \{0,1\}^{n-k-2} \\ \rho_k(w) = t}} \sk(w) \\
&= (k-2) \alpha_{k+1}(n-k-2) + \beta_{k+1}(n-k-2) \\
& \geq (k/3 - 2)\alpha_{k+1}(n-k-2).
\end{align*}

\end{proof}

In the proof of Theorem~\ref{thm:main} below, we use the following useful
inequality of Janson (see, for example, \cite{J90}).  The lower bound is standard; the upper bound is an easy modification of the one presented in \cite{AS08}.  Let $X$ be a finite set and let $P$ be a random subset of $X$, with elements $x \in X$ chosen independently with probability $p_x$.  Let $\{Z_i : i \in \cI\}$ be a system of subsets of $X$, and let $A_i$ denote the event that $Z_i \subset P$.  If $Z_i \cap Z_j = \emptyset$, then $A_i$ and $A_j$ are independent.  Let
$$
\Delta = \sum P(A_i \wedge A_j),
$$
where the sum is taken over all ordered pairs $i \neq j$ with $Z_i \cap Z_j \neq \emptyset$.  Finally, define $\mu = \sum_i P(A_i)$.

\begin{lemma} With $\mu$, $\Delta$ as above, if $\Delta \geq \mu/2$, then
$$
e^{-\mu} \leq \bigwedge_{i \in I} \overline{A_i} \leq e^{-\mu^2/3\Delta}.
$$
\end{lemma}

\begin{proof}[Proof of Theorem~\ref{thm:main}]

Suppose for the moment that $n$ is prime and $k \geq 4$.
We know that
\[
\disc(\cL_n) = \max_k (1 + \sum_{j=1}^{k-1}\sk(\ell_{n-j}) + \disc(\ell_{n-k})).\]
From Proposition~\ref{ineq}, we have that
\begin{align*}
\sum_{k=\log n + 1}^n \sk(\ell_k) &\geq \sum_{k=\log n + 1}^n (k/3-2) \cdot \alpha_{k+1}(n-k-2) \\
&\geq \sum_{k=\log n + 1}^n (k/3-2) \cdot 2^{n-k-1} (1 - n2^{-k}) \\
& = \Omega \left (\frac{2^n \log n}{n} \right ).
\end{align*}
On the other hand, for any $t$,
\begin{align*}
\sum_{k=t}^{n-1} \sk(\ell_k) &\leq \sum_{k=0}^{n-2} (2k-3) \cdot \alpha_{k+1}(n-k-2) \\
&\leq \sum_{k=1}^{n-1} 2k \cdot \alpha_k(n-k-1).
\end{align*}

We estimate this quantity using the inequality of Janson stated above.
In this case, we take $X = [n]$, $P$ is the set of indices where a $0$ appears, $p_x = 1/2$ for every $x$, $\cI = [n-k+1]$, $Z_i = [i,i+k-1]$ (i.e., the $i^\textrm{th}$ length $k$ interval of $[n]$), and $A_i$ is the event that a length $n$ word has a subsequence of the form $0^k$ on some $Z_i$.  Then
$$
\mu = (n-k+1) 2^{-k}
$$
and
\begin{align*}
\Delta &= \sum_{\substack{1 \leq i, j \leq n-k+1 \\ 0 < |i-j| < k}} 2^{-k - |i-j|} \\
& < 2^{-k+1}(n-k+1)\sum_{s=1}^\infty 2^{-s} = 2^{-k+1}(n-k+1) = 2 \mu.
\end{align*}
Furthermore,
\begin{align*}
\Delta \geq 2^{-k} (n-k+1) \sum_{s=1}^{k-1} 2^{-s} > 2^{-k-1} (n-k+1) = \mu/2,
\end{align*}
so the hypotheses hold.  Therefore, for a uniform random choice of $w \in \{0,1\}^n$,
$$
P(0^k \not \in w) \leq e^{-\mu/12} = e^{-(n-k+1)/(12 \cdot 2^{k})}.
$$
Applying this bound to the above computations,
$$
k \alpha_k(n-k-1) \leq k \cdot 2^{n-k} e^{-(n-2k)/(12 \cdot 2^{k})}.
$$
Let $T = \lfloor \log n \rfloor$.  Then
\begin{align*}
\sum_{k=1}^{n-1} k \alpha_k(n-k-1) &\leq \sum_{k=1}^{n-1} k \cdot 2^{n-k} e^{-(n-2k)/(12 \cdot 2^{k})}\\
&= 2^n \sum_{k=1}^{2 \log n} k \cdot 2^{-k} e^{-(n-2k)/(12 \cdot 2^{k})} \\
&\qquad+ 2^n \sum_{k=2 \log n+1}^{n-1} k \cdot 2^{-k} e^{-(n-2k)/(12 \cdot 2^{k})}\\
&\leq 2^n \sum_{k=-\infty}^{\infty} k \cdot 2^{-k} e^{-n/(24 \cdot 2^{k})} + o \left ( \frac{2^n \log n}{n} \right )\\
&\leq 2^n \sum_{k=-\infty}^{\infty} (T-k) \cdot 2^{k-T} e^{-n/(24 \cdot 2^{T-k})} + o \left ( \frac{2^n \log n}{n} \right )\\
&\leq 2^n \sum_{k=-\infty}^{\infty} \frac{2 \log n}{n} \cdot 2^{k} e^{-2^k/48} + o \left ( \frac{2^n \log n}{n} \right ) \\
&= O\left ( \frac{2^n \log n}{n} \right ) \cdot \sum_{k=-\infty}^{\infty} 2^{k} e^{-2^k/48} = O\left ( \frac{2^n \log n}{n} \right ).
\end{align*}
Therefore, the total discrepancy is $\Theta(2^n \log n / n)$.

There are two more terms to consider: $\sk(\ell_k)$ with $k \leq 3$, and $\max_k \disc(\ell_{n-k})$.  The former terms are bounded by $O(\rho_4^k) = O(1.93^k)$, and therefore make an insignificant contribution.  As for the latter, the length of $\ell_{n-k}$ is bounded above by $\alpha_{k+1}(n-k-2)$, and the above analysis shows that this quantity is $o(2^n \log n/n)$.  Since the length of $\ell_{n-k}$ is an upper bound for $\disc(\ell_{n-k})$, this term also does not affect the order of $\disc(\cL_n)$.

Finally, we may drop the assumption that $n$ is prime.  If not, then the above analysis is wrong: some words of length $n$, which would be part of the concatenation that gives rise to an $\ell_k$, are in fact periodic, and therefore only appear as their minimal roots in $\cL_n$.  (All Lyndon words of length dividing $n$ arise in this way.)  However, the total number of symbols they contribute is at most
$$
\sum_{d|n, d<n} d 2^{d} < n^2 2^{n/2} = o \left (\frac{2^n \log n}{n} \right ).
$$
Hence, the asymptotic bound holds.
\end{proof}

\section{Acknowledgements}

Thanks to Aaron Dutle for his careful reading of an earlier draft of this paper, and to Ron Graham for suggesting the problem.


\begin{thebibliography}{100}
\bibitem{AS08} N.~Alon, J.~Spencer, The probabilistic method. Wiley-Interscience Series in Discrete Mathematics and Optimization. John Wiley \& Sons, Inc., Hoboken, NJ, 2008.
\bibitem{F57} L.~R.~Ford, A cyclic arrangement of $m$-tuples, Report P-1071, Rand Corp., Santa Monica, CA, 1957.
\bibitem{F82} H.~Fredricksen, A survey of full length nonlinear shift register cycle algorithms, {\it SIAM Rev.}~{\bf 24} (1982), no.~2, 195--221.
\bibitem{H90} Y.~J.~Huang, A new algorithm for the generation of binary de Bruijn sequences, {\it J.~Algorithms} {\bf 11} (1990), no.~1, 44--51.
\bibitem{J90} S.~Janson, Poisson approximation for large deviations, {\it Random Struct. Alg.} {\bf 1}, 221--229 (1990).
\bibitem{L01} J.~H.~van Lint, R.~M.~Wilson, A course in combinatorics. Second edition. Cambridge University Press, Cambridge, 2001.
\bibitem{M34} M.~H.~Martin, A problem in arrangements, {\it Bull.~Amer.~Math.~Soc.} {\bf 40} (1934), 859--864.
\bibitem{O95} A.~M.~Odlyzko, Asymptotic enumeration methods. {\it Handbook of combinatorics}, Vol. 2, 1063--1229, Elsevier, Amsterdam, 1995.
\end{thebibliography}
\end{document}